\documentclass[11pt,a4paper]{article}
\usepackage{graphicx}
\usepackage{verbatim}
\usepackage{amsthm}
\usepackage{amsmath}
\usepackage{setspace}
\usepackage{amssymb}
\newtheorem{theorem}{Theorem}

\newtheorem{lemma}{Lemma}

\newcommand{\poi}{ Poisson process with a negative exponential component}

\begin{document}
\title
{On several  two-boundary problems  for a  particular  class of
L\'{e}vy processes}
\date{}
\author {T. Kadankova  \thanks{
    Hasselt University,  Center for Statistics,  Agoralaan, Building D,
3590 Diepenbeek, Belgium, \newline tel.: +32(0)11 26 82 97,\:
e-mail: tetyana.kadankova@uhasselt.be} \:and \: N. Veraverbeke
\\Hasselt University,  Belgium} \maketitle


\noindent{\bf Running head:} Two-boundary problems for certain
L\'{e}vy processes

\begin{abstract}
Several two-boundary problems   are  solved  for a special  L\'{e}vy
process:   the Poisson process with an exponential component. The
jumps of this process  are controlled   by a homogeneous Poisson
process, the positive   jump size  distribution is arbitrary, while
the distribution of the negative jumps is  exponential. Closed  form
expressions  are obtained  for  the integral transforms of the joint
distribution of the first exit time  from an interval and the value
of the overshoot through boundaries at the first exit time. Also the
joint distribution of the  first entry time  into the interval and
the value of the process at this time instant are determined in
terms of integral transforms.
\end{abstract}
\section{ Introduction}
We assume that all random variables and  stochastic processes  are
defined  on   $(\Omega, \: \mathfrak{F}, \{\mathfrak{F}_t\}, P),$ a
filtered probability space, where the filtration $
\{\mathfrak{F}_t\}$ satisfies  the usual conditions of
right-continuity and completion. A L\'{e}vy process is  a
$\mathfrak{F}$-adapted stochastic process $\;\{ \xi(t);\: t\ge0 \}$
which has independent and stationary increments and whose  paths are
right-continuous  with left limits \cite{Bert1996}. Under the
assumption that $\xi(0)=0$ the Laplace transform of the process
$\;\{ \xi(t);\: t\ge0 \}$ has the form
$E[e^{-p\xi(t)}]=e^{t\:k(p)},$ $ {\rm Re}\,p=0,$ where the function
$k(p)  $ is called  the Laplace exponent and is given by the formula
(\cite{Sk2})
\begin{align}
k(p)=\frac{1}{t}\ln E[ e^{-p\xi(t)}]= \frac1{2}\,p^2\sigma^2-\alpha
p+\int_{-\infty}^{\infty}
\left(e^{-px}-1+\frac{px}{1+x^2}\right)\Pi(dx).
\end{align}
Here  $ \alpha,\: \sigma \in \mathbb{R}$ and $\Pi(\cdot) $ is a
measure on  the real line. The introduced process is a space
homogeneous,  strong Markov process. Note, that  the distribution of
the first exit time from  an interval  plays a crucial role in
applications and its  knowledge  also allows to solve a number of
other two-boundary problems. Let us fix $B>0$ and define the
variable
$$
 \chi(y)=\inf\{\,t:\;y+\xi(t)\notin[0,B]\,\},\quad y\in [0,B], \:
$$
the first exit time from the interval $[0,B]$ by the process
$y+\xi(t).$ The random variable $\,\chi(y)\,$  is a Markov time and
$P\,[\,\chi(y)<\infty\,]=1$  \cite{Sk2}. Exit from the interval
$[0,B]$ can take place either through the upper boundary  $ B, $ or
through the lower boundary  $0.$ Introduce events:
$A^{\,B}=\{\omega:\,\xi(\chi(y))> B\,\}, $ i.e.  the exit takes
place through the upper boundary; $
A_{0}=\{\omega:\,\xi(\chi(y))<0\,\},  $  i. e. the exit takes place
through the lower boundary. Define
$$
X(y)=(\xi(\chi(y))-B)\, I_{\,A^{\,B}}+(-\xi(\chi(y)))\, I_{\,A_{0}},
\qquad  P\,[\,A^{\,B}+A_{0}\,]=1,
$$
the value of the overshoot through  one of the boundaries at  the
first  exit time, where  $ I_{\,A}= I_{\,A}(\omega) $ is the
indicator of the event  $ A.$ The first two-boundary problem  for
L\'{e}vy processes with the Laplace exponent of the general form (1)
has been solved by Gihman and   Skorokhod (\cite{Sk2}, p.306-311).
These authors have determined the joint distribution of
$\{\,\xi^{-}(t),\xi(t),\xi^{+}(t)\,\},$ where $ \xi^{+}(t)=\sup_{u
\le t}\limits\,\xi(u),$ $ \xi^{-}(t)=\inf_{u \le t}\limits\,\xi(u),
\: t \ge 0. $ For a spectrally positive L\'{e}vy process the joint
distribution of  $\{\chi(y),X(y)\,\}$  has been studied by many
authors among which  Emery \cite{Emery},  Suprun
 and Shurenkov  \cite{SuSh}. The first exit time  for a spectrally one-sided
L\'{e}vy process has been considered by Bertoin \cite{Bert1997},
Pistorius \cite{Pistor2004}, \cite{Pist2004}, Kyprianou
\cite{Kypr2003} and others. Kadankov and Kadankova \cite{2Ka2} have
suggested another approach for determining the joint distribution of
$\{\chi(y),X(y)\}$  for the L\'{e}vy process with Laplace exponent
(1). Their method is based on application of one-boundary
functionals  $ \{\tau^{\,x},\,T^{\,x}\},$ $ \{\tau_{x},\,T_{x}\},$
$x\ge0,$ where
$$
\tau^{\,x}=\inf\{t:\,\xi(t)> x\},\; T^{\,x}=\xi(\tau^x)-x,\quad
\tau_x=\inf\{t:\,\xi(t)< -x\},\; T_x=-\xi(\tau_x)-x.
$$
Integral transforms of these joint distributions have been obtained
in 60's  in papers of  Pecherskii and   Rogozin
\cite{Rog1},\cite{PeRog}, Borovkov \cite{Br1},  Zolotarev
\cite{Zol1}. Kadankov and Kadankova \cite{2Ka2} have used
probabilistic methods (the total probability law, space homogeneity
and  the strong Markov property of the process) to determine   the
integral transforms $ E[e^{-s\chi(y)};X(y)\in du,A^{\,B}],$ $
E[e^{-s\chi(y)};X(y)\in du,A_{0}] $ of the joint distribution of
$\{\chi(y),X(y)\}.$  For a spectrally positive L\'{e}vy  process
several two-boundary problems have been solved in \cite{
Kad6}-\cite{Kad2}.

In this paper  we obtain  the integral  transforms  of the
distributions  of a number  of two-boundary  functionals associated
with the first  exit time  (Section 3)  and the first entry time
(Section 4)  for an important  particular class of L\'{e}vy
processes described in Section 2. The advantage   is that these are
closed formulas for the transforms, with no recursions, typical for
the general case.
\section{The Poisson process  with  negative  exponential exponent}
We now  give a formal definition of the  process which we consider.
Let $\eta \in(0,\infty)$  be a positive random variable, and   
$\gamma $  be an exponential variable with  parameter $\lambda>0:$ $
P\,[\,\gamma>x\,]=e^{-\lambda x},$ $x\ge0.$ Introduce the random
variable  $\xi\in \mathbb{R} $ by its distribution function
$$
F(x)=a\,e^{x\lambda}\, I\{x\le0\}+(a+(1-a)\, P\,[\,\eta \le x\,])\,
I\{x>0\}, \qquad a\in(0,1),\quad \lambda>0.
$$
Consider a right-continuous compound Poisson process $
\xi(t)=\sum_{k=0}^{N(t)}\: \xi_k,$ $  t \ge 0, $  where $\{ \xi_k;\:
k \ge 1\} $ are independent random variables identically distributed
as  $\xi,$  $\xi_0=0,$ and $ N(t)$ is a homogeneous Poisson process
with intensity $ c>0. $ Then its Laplace exponent is of the form
\begin{align} k(p)=c\int_{-\infty}^{\infty}(e^{-xp}-1)\,dF(x)=
a_1\frac{p}{\lambda-p}+a_2( E[\,e^{-p\eta}]-1),                        
\qquad c>0,\quad {\rm Re}\,p=0,
\end{align}
where $a_1=ac,$ $a_2=(1-a)c.$  Here and in the sequel we will call
such process  the Poisson process  with a  negative exponential
component. Note, that  inter-arrival times of the jumps of the
process $\{\xi(t);\: t\ge0\}$ are exponentially distributed with
parameter  $c.$  With probability $1-a $ there occur positive jumps
of size  $\eta,$ and with probability $a$ there occur negative jumps
of value  $\gamma$ that is exponentially distributed with parameter
$\lambda.$ The first term of (2) is the simplest case of a rational
function,  while the second term is nothing but the Laplace exponent
of a monotone Poisson process with positive jumps of value $\eta.$
It is well known fact (\cite{BGu}), that in this case the equation
$k(p)-s=0,$ $s>0 $ has a unique root $ c(s)\in(0,\lambda),$ in the
semi-plane $ {\rm Re}\,p>0. $ Denote by $\nu_s$ an exponentially
distributed random variable with parameter $s>0,$  independent of
the process. The introduction of  $\nu_s$ allows  the following
short-hand notation  for the  double Laplace transforms  of the
process, i.e. $\int\limits_{0}^{\infty}se^{-st}E[e^{-p
\xi(t)}]dt=E[e^{-p \xi(\nu_s)}]. $ For the double  Laplace
transforms of the processes $\xi^{+}(\cdot),$ $\xi^{-}(\cdot)$  the
following formulae hold
\begin{align}
&
E[\,e^{-p\xi^{-}(\nu_s)}]=\frac{c(s)}{\lambda}\,\frac{\lambda-p}{c(s)-p},
\qquad {\rm Re}\,p\le0,\notag\\
& E[\,e^{-p\xi^{+}(\nu_s)}]=\frac{s\lambda}{c(s)}\,(p-c(s))\,\mathbb R(p,s),  
\qquad {\rm Re}\,p\ge0,
\end{align}
where
\begin{align}
\mathbb R(p,s)=\left(a_1p+(p-\lambda)
[s-a_2(\,E[\,e^{-p\eta}]-1)]\right)^{-1},                               
\qquad {\rm Re}\,p\ge0, \quad p\ne c(s).
\end{align}
Observe that the  function  $\mathbb R(p,s) $   is analytic in the
semi-plane $ {\rm Re}\,p>c(s)$ and $\lim_{p\to\infty}\limits\mathbb
R(p,s)=0.$ Therefore,  it allows a representation in the  form of an
absolutely convergent Laplace integral  (\cite{DiP1})
\begin{align}
\mathbb R(p,s)=\int_{0}^{\infty}e^{-px}R_x(s)\,dx,
\qquad {\rm Re}\,p>c(s).                                      
\end{align}
We will call the function   $ R_x(s),$ $x\ge0$  the resolvent of the
Poisson process  with a  negative exponential  component. We assume
that  $R_x(s)=0,$ for $x<0.$ Note, that   $
R_0(s)=\lim_{p\to\infty}\limits p\,\mathbb R(p,s)=(c+s)^{-1}, $ and
the equalities (3)  imply
\begin{align*}
 P\,[\,\xi^{-}(\nu_s)=0\,]=\frac{c(s)}{\lambda},\qquad
P\,[\,\xi^{+}(\nu_s)=0\,]=\frac{\lambda}{c(s)}\,\frac{s}{s+c}.
\end{align*}
The second formula of  (3) yields
\begin{align}
\mathbb R(p,s)= \frac{c(s)}{s\lambda}\,\frac{1}{p-c(s)}\,
E[\,e^{-p\xi^{+}(\nu_s)}],
\qquad {\rm Re}\,p>c(s).                                 
\end{align}
The functions
\begin{align*}
\frac{1}{p-c(s)}&=\int_{0}^{\infty}e^{-u(p-c(s))}du,\qquad {\rm
Re}\,p>c(s), \quad \\
E[\,e^{-p\xi^{+}(\nu_s)}]&=\int_{0}^{\infty}e^{-up} d\,
P\,[\,\xi^+(\nu_s)<u\,],\qquad {\rm Re}\,p\ge0,
\end{align*}
which enter the right-hand side of (6),  are the Laplace transforms
for ${\rm Re}\,p>c(s).$ Therefore,  the original functions  of the
left-hand side and the right-hand side of (6) coincide, and
\begin{align}
R_x(s)= \frac{c(s)}{s\lambda}\,\int_{-0}^{x}e^{c(s)(x-u)}
d\, P\,[\,\xi^+(\nu_s)<u\,],\qquad x\ge0,                     
\end{align}
which is  the resolvent  representation   of  the Poisson process
with a  negative exponential  component. Note, that the
representation for the resolvent  of the  spectrally one-sided
L\'{e}vy process similar to (7) was obtained  by Suprun and
Shurenkov \cite{Su},\cite{SuSh}. This representation implies that
$R_x(s),$ $x\ge0$ is a positive, monotone,  continuous, increasing
function of an exponential order, i.e. there exists $ 0 <
A(s)<\infty$ such that $ R_x(s)<A(s)e^{xc(s)},$ for all $x\ge0.$
Therefore,
$$
\int_{0}^{\infty}R_x(s)e^{-\alpha x}dx<\infty,\qquad \alpha>c(s).
$$
Moreover,  in  the neighborhood  of  any $x\ge0 $ the function
$R_x(s)$ has bounded variation. Hence, the inversion formula  is
valid
\begin{align}
R_x(s)=\frac{1}{2\pi
i}\int_{\alpha-i\infty}^{\alpha+i\infty}e^{xp}\,\mathbb R(p,s)\,dp,
\qquad \alpha>c(s).                                              
\end{align}
The latter equality together with  (5) determines  the resolvent of
the Poisson process with a  negative exponential  component. To
derive the joint distribution  of the first exit time and the value
of the overshoot at the first  exit time  for a \poi, \: we apply a
general theorem for the  L\'{e}vy processes which has been proved in
\cite{2Ka2}. Before stating the theorem we mention the following
results 
\begin{align}
& E\,[\,e^{\,-s\tau^x-pT^x}\,]= \left(                                   
E\,[e^{-p\xi^{+}(\nu_s)}]\right)^{-1}
E\,[\,e^{-p(\xi^{+}(\nu_s)-x)}; \; \xi^{+}(\nu_s)> x\,],
& \qquad {\rm Re}\,p \ge 0,\notag\\
& E\,[\,e^{\,-s\tau_x -pT_x}\,]= \left(
E\,[e^{\,p\xi^{-}(\nu_s)}]\right)^{-1}
E\,[\,e^{p(\xi^{-}(\nu_s)+x)}; \;-\xi^{-}(\nu_s)> x\,], &\qquad {\rm
Re}\,p \ge 0.
\end{align}
The formulae (9) have been obtained  by  Pecherskii  and Rogozin
\cite{PeRog}. A simple proof of these equalities is given in
\cite{2Ka2}. After some calculations it follows from  (3)  and (9)
 that the integral transforms of the joint
distributions of  $\{\tau_x,\,T_x\},$ $\{\tau^x,\,T^x\}$ of the
Poisson process with a negative exponential component satisfy the
equalities
\begin{align}
& E[\,e^{-s\tau_{x}}; \, T_{x}\in du\,]=(\lambda-c(s))\,e^{-xc(s)}\,
e^{-\lambda u}\,du= E[\,e^{-s\tau_{x}}]\, P\,[\,\gamma\in du\,],\\
&\int_{0}^{\infty}e^{-px}\,                                        
 E [\,e^{\,-s\tau ^{x}-z\xi(\tau^{x})}]\,dx=
\frac{1}{p}\left(1-\frac{p+z-c(s)}{z-c(s)}\; \frac{\mathbb
R(p+z,s)}{\mathbb R(z,s)}\right), \quad {\rm Re}\,p>0,\;{\rm
Re}\,z\ge0.\notag
\end{align}
The first  equality  of (10) yields that   $\tau_x$ and  $T_x$ are
independent. Moreover,  for all  $x\ge0 $  the  value of the
overshoot  through the lower level   $T_x $ is  exponentially
distributed  with  parameter $\lambda.$ This  feature characterizes
the Poisson process with a negative exponential  component.  Now we
state the main results on two-sided exit problems.
\section{  The first exit from  an  interval}
We now derive the joint distribution of the first exit time and the
value of the overshoot at the first exit.  The following result  is
true for the  general L\'{e}vy processes (\cite{2Ka2}),  and for
convenience  it is stated   here as a lemma.
\begin{lemma}                           
Let    $\{\xi(t);\: t \ge 0\} ,$ $ \xi(0)=0$ be a  real-valued
L\'{e}vy process whose  Laplace exponent  is given by (1), $B>0,$
$\,y\in[0,B],\quad x=B-y.$ Let
$$
 \chi(y)=\inf\{\,t:\,y+\xi(t)\notin[0,B]\,\},\qquad
X(y)=(\xi(\chi(y))-B)\, I_{\,A^{\,B}}+(-\xi(\chi(y)))\,
I_{\,A_{\,0}}
$$
be, respectively, the first exit  time   from the interval $[0,B]$
by the process $y+\xi(t)$ and the value of the overshoot through the
boundary at the first  exit time.  For $s>0,$ the Laplace transforms
of the joint distribution of $\{\,\chi(y),X(y)\,\} $ satisfy the
following equations
\begin{align}
& E\,[e^{-s\chi(y)};X(y)\in du,A^{\,B}]= f_{+}^{s}(x,du)+
\int_{0}^{\infty}f_{+}^{s}(x,dv)\, K_{+}^{s}(v,du),\notag\\
& E\,[e^{-s\chi(y)};X(y)\in du,A_{\,0}]= f_{-}^{s}(y,du)+
\int_{0}^{\infty}f_{-}^{s}(y,dv)\, K_{-}^{s}(v,du),                           
\end{align}
where
\begin{align*}
& f_{+}^{s}(x,du)=  E\,[e^{-s\tau^x};T^x\in du] -\int_{0}^{\infty}\,
E\,[e^{-s\tau_y};T_y\in dv]
\,  E\,[\,e^{-s\tau^{v+B}};T^{v+B}\in du],\\
& f_{-}^{s}(y,du)=  E\,[e^{-s\tau_y};T_y\in du] -\int_{0}^{\infty}\,
E\,[e^{-s\tau^x};T^x\in dv] \, E\,[e^{-s\tau_{v+B}};T_{v+B}\in du];
\end{align*}
and $
K_{\pm}^s(v,du)=\sum_{n=1}^{\infty}\limits\,K_{\pm}^{(n)}(v,du,s),
\quad v\ge0 $ are the series of the successive iterations of kernels
 $K_{\pm}(v, du,s).$ These kernels are given by
 \begin{align}
&K_+(v,du,s)=\int_{0}^{\infty} E[e^{-s\tau_{v+B}};T_{v+B}\in dl]
\, E[e^{-s\tau^{l+B}};T^{l+B}\in du], \notag\\
&K_-(v,du,s)=\int_{0}^{\infty} E[e^{-s\tau^{v+B}};T^{v+B}\in          
dl] \, E[e^{-s\tau_{l+B}};T_{l+B}\in du],
 \end{align}
and  their successive iterations $(n\in \mathbb{N}=\{1,2, \dots\})$
are defined by
\begin{align}
K_{\pm}^{(1)}(v,du,s)=K_{\pm}(v,du,s),\quad K_{\pm}^{(n+1)}(v,du,s)=
\int_{0}^{\infty} K_{\pm}^{(n)}(v,dl,s)\,K_{\pm}(l,du,s).      
\end{align}
\end{lemma}
We apply now the formulae  of Lemma 1  for the case  when the
underlying process is  the Poisson process with  an exponentially
distributed negative component.
\begin{theorem}                                   
Let   $ \{\xi(t);\:t\ge0\}, \: \xi(0)=0 $ be a  Poisson process with
a negative exponential  component whose  Laplace exponent  is given
by (2), $ B> 0,$ $ y\in[0,B],$ $ x=B-y.$  Let
$$
 \chi(y)=\inf\{\,t:\,y+\xi(t)\notin [0,B]\,\},\qquad
  X(y)=(\xi(\chi(y))-B)\, I_{\,A^{\,B}}+
  (-\xi(\chi(y)))\, I_{\,A_{\,0}}
$$
be, respectively, the  first exit time  from the interval and the
value  of overshoot through   one of the boundaries. Then  for $
s>0, $

1) the integral transforms of the joint distribution of
 $\{\,\chi(y),\;X(y)\,\} $  satisfy the following equations
\begin{align}
& E[e^{-s\chi(y)};X(y)\in du,A_{0}]= e^{-\lambda u}
(\lambda-c(s))\,e^{-yc(s)}
\left(1- E\,[e^{-s\tau^x-c(s)\xi(\tau^x)}]\right)K(s)^{-1}du,\\
& E[e^{-s\chi(y)};X(y)\in du,A^{\,B}]=
E[e^{-s\tau^x};T^x\in du] - E[e^{-s\chi(y)};A_{0}]\,                  
E[e^{-s\tau^{\gamma+B}};T^{\gamma+B}\in du],\notag
\end{align}
where
\begin{align*}
&K(s)=1- E[e^{-s\tau_B}]\, E[e^{-s\tau^{\gamma+B}-c(s)T^{\gamma+B}}]\\
& E[\,e^{-s\tau^{\gamma+B}-c(s)T^{\gamma+B}}\,]=
\lambda\int_{0}^{\infty}e^{-\lambda u}
E\,[\,e^{-s\tau^{u+B}-c(s)T^{u+B}}\,]\,du.
\end{align*}
In particular,
\begin{align}                                                       
 E[e^{-s\chi(y)};A_{\,0}]=
&\left(1-\frac{c(s)}{\lambda}\right)\,e^{-yc(s)}
\left(1- E[e^{-s\tau^x-c(s)\xi(\tau^x)}]\right)K(s)^{-1}, \\
 E[e^{-s\chi(y)};A^{\,B}]= & E[\,e^{-s\tau^x}]-
E[e^{-s\chi(y)};A_{\,0}]\, E[\,e^{-s\tau^{\gamma+B}}];\notag
\end{align}

2)  the Laplace transforms of the random variable  $ \chi(y) $
satisfy  the following representations
\begin{align}
 E[e^{-s\chi(y)};X(y)\in du,A_{0}] =&
e^{-\lambda(u+B)}\, \frac{R_{x}(s)}{\hat
R_{B}(\lambda,s)}\,du,\qquad
E[e^{-s\chi(y)};A_{0}]=\frac1{\lambda}\,e^{-\lambda B}\,
\frac{R_{x}(s)}{\hat R_{B}(\lambda,s)},\notag\\
 E[e^{-s\chi(y)};A^{\,B}] =&1-\frac{R_{x}(s)}{\hat R_B
(\lambda,s)}\, \left[\frac{1}{\lambda}\,e^{-\lambda B}
+s\lambda \,\hat S_B (\lambda,s) \right]+s\lambda\,S_{x}(s),\notag\\
\int_{0}^{\infty}e^{-st} P[\chi(y)>t]\,dt=&
\lambda\,\frac{R_{x}(s)}{\hat R_B (\lambda,s)}\, \hat S_B                                  
(\lambda,s) -\lambda\,S_{x}(s),
\end{align}
where $ R_x(s),$ $x\ge0$   is the resolvent  of the process, defined
by   (5), (8);
$$
S_x (s)=\int_{0}^{x}R_u (s)\,du, \qquad \hat
R_{B}(\lambda,s)=\int_{B}^{\infty}e^{-\lambda u}  R_u (s)\,du,\qquad
\hat S_{B}(\lambda,s) =\int_{B}^{\infty}e^{-\lambda u} S_u (s)\,du.
$$
\end{theorem}
\begin{proof}
\rm For the Poisson process with a negative exponential component,
equalities of Lemma  1 take a simplified form. Using the equalities
(10) and the defining formulae  (12)  for the kernels
$K_{\pm}(v,du,s)$  we obtain
\begin{align*}
&K_{+}(v,du,s)
=\left(1-\frac{c(s)}{\lambda}\right)\,e^{-c(s)(v+B)}\,
 E[e^{-s\tau^{\gamma+B}};T^{\gamma+B}\in du],\\
&K_{-}(v,du,s) =e^{-\lambda u}(\lambda-c(s))\,e^{-c(s)B}\,
E[e^{-s\tau^{v+B}-c(s)T^{v+B}}]\,du,
\end{align*}
where $\gamma $ is an exponentially distributed random variable with
the parameter  $\lambda,$ independent of the  process under  the
consideration. Using these equalities, the method of mathematical
induction  and the formulae (13), we obtain the successive
iterations $K_{\pm}^{(n)}(v,du,s),$ $n\in\mathbb{N} $ of the kernels
$K_{\pm}(v,du,s):$
\begin{align*}
K_{-}^{(n)}(v,du,s)=& E[e^{-s\tau^{v+B}-c(s)T^{v+B}}] \left(
E\,[e^{-s\tau_{B}}] \right)^n \left(
E[e^{-s\tau^{\gamma+B}-c(s)T^{\gamma+B}}]\right)^{n-1}\;
\lambda e^{-\lambda u}\,du,\\
K_+^{(n)}(v,du,s)=&e^{-vc(s)}\left( E\,[e^{-s\tau_{B}}] \right)^n
\left( E[\,e^{-s\tau^{\gamma+B}-c(s)T^{\gamma+B}}\,]\right)^{n-1}
E[e^{-s\tau^{\gamma+B}};T^{\gamma+B}\in du\,].
\end{align*}
The series  $ K_{\pm}^{s}(v,du) $ of the successive iterations
$K_{\pm}^{\,(n)}(v,du,s)$  are nothing but geometric series,  and
their sums  are given by
\begin{align*}
& K_-^{s}(v,du)=\sum_{n=1}^{\infty}K_{-}^{(n)}(v,du,s)=
E[e^{-s\tau^{v+B}-c(s)T^{v+B}}]\, E[e^{-s\tau_{B}}]
\;K(s)^{-1}\lambda\:e^{-\lambda u}\,du,\\
& K_+^{s}(v,du)=\sum_{n=1}^{\infty}K_{+}^{(n)}(v,du,s)= e^{-vc(s)}\,
E\,[e^{-s\tau_{B}}]\; K(s)^{-1}\;
E[e^{-s\tau^{\gamma+B}};T^{\gamma+B}\in du].
\end{align*}
Substituting the obtained  expressions for  $ {K}_{\pm}^{s}(v,du)$
into  (11) yields  the formulae  (14) of Theorem 1. Integrating the
formulae (14) with respect to  $u\in \mathbb{R}_+$  leads to formula
(15) of the theorem. Now, utilizing the definition of the resolvent
(5), (8) and the equalities (10)  we derive the resolvent
representation for the functions $
E[e^{-s\tau^{\,x}-c(s)\xi(\tau^{x})}],$ $ E[e^{-s\tau^x}]:$
\begin{align*}
& E[e^{-s\tau^{\,x}-c(s)\xi(\tau^{x})}]=
1-e^{-xc(s)}\;R_x(s)\;r(c(s),s),\\
& E[e^{-s\tau^x}]=1-\frac{s\lambda }{c(s)}\,R_x(s)+
s\lambda\,S_x(s),
\end{align*}
where
$$
S_x(s)=\int_{0}^x\,\,R_u(s)\,du,\qquad
r(c(s),s)=\left.\frac{d}{dp}\mathbb R(p,s)^{-1}\right|_{p=c(s)}.
$$
Substituting  these  expressions  into (15), we  obtain the
representations (16) of Theorem 1.
\end{proof}
\section{ The first entry time into an interval}
The knowledge  of  the joint distribution  of $\{\chi(y),\, X(y)\}$
allows  to solve another two-boundary problem, namely  to determine
the integral transforms of the joint distribution of the first entry
time  into the fixed interval by the L\'{e}vy process and the value
of the process at this time. We  obtain the result   in Theorem 2
below for the  general L\'{e}vy processes. In Theorem 3 closed
expressions   for the  integral transforms  in case of a Poisson
process with a negative exponential exponent are given.
\begin{theorem}                                         
Let $\{\xi(t);\: t\ge0\},$ $\xi(0)=0$ be a  L\'{e}vy process whose
Laplace exponent is given by  (1), $B>0,$ $\chi(y)\stackrel{{\rm
def}}{=}0$ for $y\notin [0,B].$  Let
$$
\overline\chi(y)=\inf\{t>\chi(y):\,y+\xi(t)\in [0,B]\},\qquad
\overline X(y)=y+\xi(\overline\chi(y))\in [0,B],\quad y\in
\mathbb{R}
$$
be, respectively,  the   first  entry time   into the interval $
[0,B] $ by the process $ y+\xi(t) $  and  the value of the process
at this time.  For $s>0,$ the integral transforms of the joint
distribution of $\{\overline\chi(y),\overline X(y)\},$ $y\in
\mathbb{R} $ satisfy the following equations
\begin{align}
b^v(du,s)& \stackrel{def}{=} E[e^{-s\overline\chi(v+B)};\overline
X(v+B)\in du]= \int_{0}^{\infty} Q_+^s(v,dl)\,
 E[e^{-s\tau_{l}};B-T_{l}\in du] \notag\\
&+\int_{0}^{\infty} Q_+^s(v,dl)\,\int_{0}^{\infty}                     
 E[e^{-s\tau_{l}};T_{l}-B\in d\nu] \,
 E[e^{-s\tau^{\nu}};T^{\nu}\in du],\qquad  v>0,\notag\\
b_v(du,s)&\stackrel{def}{=} E[e^{-s\overline\chi(-v)};\overline
X(-v)\in du]= \int_{0}^{\infty} Q_-^s(v,dl)\,
 E[e^{-s\tau^{l}};T^{l}\in du] \\
&+\int_{0}^{\infty} Q_-^s(v,dl)\,\int_{0}^{\infty}
E[e^{-s\tau^{l}};T^{l}-B\in d\nu] \,
 E[e^{-s\tau_{\nu}};B-T_{\nu}\in du],\qquad  v>0,\notag\\
b(y,du,s)&\stackrel{def}{=} E[e^{-s\overline\chi(y)};\overline
X(y)\in du]
=\int_{0}^{\infty} E[e^{-s\chi(y)};X(y)\in dv,A^{\,B}]\,b^v(du,s)\notag\\
&+\int_{0}^{\infty} E[e^{-s\chi(y)};X(y)\in dv,A_{\,0}]\,b_v(du,s),
\qquad  y\in [0,B],\notag
\end{align}
where $ \delta(x),$ $x\in \mathbb{R}$ is the delta function and
\begin{align}
 Q_{\pm}^s(v,du)=\delta(v-u)\,du+\sum\limits_{n=1}^{\infty}               
\limits Q_{\pm}^{(n)}(v,du,s),\qquad v>0.
\end{align}
The functions   $Q_{\pm}^{(n)}(v,du,s),$ $n\in\mathbb{N}$   are
defined  by
\begin{align}
Q_{\pm}^{(1)}(v,du,s)=Q_{\pm}(v,du,s),\qquad
Q_{\pm}^{(n+1)}(v,du,s)=\int_{0}^{\infty}
Q_{\pm}^{(n)}(v,dl,s)Q_{\pm}(l,du,s);                        
\end{align}
and they  are the successive iterations of the kernels
$Q_{\pm}(v,du,s),$ which are given by
\begin{align}
&Q_+(v,du,s)=\int_{0}^{\infty}  E[e^{-s\tau_{v}};T_{v}-B\in dl]  
\, E[e^{-s\tau^{l}};T^{l}-B\in du],\notag\\
&Q_-(v,du,s)=\int_{0}^{\infty}  E[e^{-s\tau^{v}};T^{v}-B\in dl] \,
E[e^{-s\tau_{l}};T_{l}-B\in du].
\end{align}
\end{theorem}
\begin{proof}\rm
For the functions $b^v(du,s),$ $b_v(du,s),$ $v>0$ according to the
total probability  law, space homogeneity  of the process  and the
fact that  $\tau_v,$ $\tau^v$ are Markov times, the following system
of equations  is valid
\begin{align}
b^v(du,s)=& E[e^{-s\tau_v};B-T_v\in du]+
\int_{0}^{\infty} E[e^{-s\tau_v};T_v-B\in dl]\,b_{l}(du,s),\notag\\
b_v(du,s)=& E[e^{-s\tau^v};T^v\in du]+
\int_{0}^{\infty} E[e^{-s\tau^v};T^v-B\in dl]\,b^{l}(du,s).    
\end{align}
This system is similar to a system of linear equations  with two
variables. Substituting the expression for   $b_v(du,s)$ from the
right-hand side of the second equation into the first equation
 yields
 \begin{align*}
b^v(du,s)&= E[e^{-s\tau_v};B-T_v\in du]+ \int_{0}^{\infty}
E\,[e^{-s\tau_v};T_v-B\in dl]\,
 E[e^{-s\tau^l};T^l\in du]\\
&+\int_{l=0}^{\infty} E[e^{-s\tau_v};T_v-B\in dl]\,
\int_{\nu=0}^{\infty} E[e^{-s\tau^l};T^l-B\in d\nu]\,b^\nu(du,s).
\end{align*}
Changing the order of  integration  in the third term of the second
equation  implies for the function $ b^v(du,s),$ $v>0$
\begin{align}
b^v(du,s)&=\int_{0}^{\infty}Q_{+}(v,d\nu,s)\,b^\nu(du,s)\\
&+ E\,[\,e^{-s\tau_v};B-T_v\in du\,]+
\int_{0}^{\infty} E[e^{-s\tau_v};T_v-B\in dl]\,                         
 E[e^{-s\tau^l};T^l\in du],\notag
\end{align}
which is a linear integral  equation  with the following kernel
\begin{align*}
Q_+(v,du,s)=\int_{0}^{\infty}  E[e^{-s\tau_{v}};T_{v}-B\in dl] \,
E[e^{-s\tau^{l}};T^{l}-B\in du],\qquad v>0.
\end{align*}
We now show, that for all $v,u>0,$ $s>s_0>0 $ this kernel enjoys the
following property
$$
Q_+(v,du,s)<\lambda,\qquad \lambda= E[\,e^{-s\tau_B}]\,
E[\,e^{-s\tau^B}], \quad s_0>0.
$$
Indeed, for all $s>0$  it follows from
\begin{align*}
 E[e^{-s\tau^v};T^v  -B\in du]=
 E[e^{-s\tau^{v+B}};T^{v+B}\in du]
-\int_{0}^{B}E\,[e^{-s\tau^v};T^v\in dl]\,
E[e^{-s\tau^{B-l}};T^{B-l}\in du],
\end{align*}
that the following chain of inequalities holds
\begin{align*}
E[e^{-s\tau^v};T^v-B\in du] \le E[e^{-s\tau^{v+B}};T^{v+B}\in du]
\le E[e^{-s\tau^{v+B}}] \le E[e^{-s\tau^B}].
\end{align*}
Analogously we establish, that $ E[e^{-s\tau_v};T_v-B\in du] \le
E[e^{-s\tau_{v+B}};T_{v+B}\in du] \le E[e^{-s\tau_{v+B}}] \le
E[e^{-s\tau_B}].$
 These chains of the inequalities imply the
following estimation for the kernel $Q_+(v,du,s),$ for all $v,u>0,$
$s>s_0>0 $
\begin{align*}
Q_+(v,du,s)&=\int_{0}^{\infty}  E[e^{-s\tau_{v}};T_{v}-B\in dl]
\, E[e^{-s\tau^{l}};T^{l}-B\in du]\\
&\le E[e^{-s\tau_B}]\, E[e^{-s\tau^B}]<\lambda= E[e^{-s_0\tau_B}]\,
E[e^{-s_0\tau^B}],\qquad s_0>0.
\end{align*}
This bound and the method of mathematical induction  yield that the
successive iterations $ Q_+^{\,(n)}(v,du,s) $  (19) of the kernels  
$ Q_+(v,du,s)$ for all  $v,u>0,$ $s>s_0>0 $ obey the inequality
$Q_{+}^{(n+1)}(v,du,s) < \lambda^{n+1},\: n\in\mathbb{N}. $
Therefore, the series of successive iterations $ \sum_{n=1}^{\infty}
\limits Q_{+}^{(n)}(v,du,s)< \lambda(1-\lambda)^{-1} $ converges
uniformly for all  $v,u>0,$ $s>s_0>0.$ Utilizing now the method of
successive iterations (\cite{Pet1}, p. 33) to solve the integral
equation (22) yields the first equality of the theorem. The second
equality  of the theorem can be verified analogously. It is not
difficult to establish the third equality of the theorem using  the
total probability law and the fact that $\chi(y)$ is  the Markov
time.
\end{proof}
 Denote by
\begin{align*}
m_{\gamma}^s(du)=\int_0^{\infty}\lambda\,e^{-\lambda x}
E[\,e^{-s\tau^x};\;T^x\in du\,]\,dx,\qquad P(\lambda,du)=e^{-\lambda
B}(m_{\gamma}^s(du)+\lambda\, e^{\lambda u}du).
\end{align*}
\begin{theorem}                                         
Let  $\{\xi(t);\: t\ge0\},$ $\xi(0)=0$  be a Poisson process with a
negative exponential  component as specified  above, $B>0,$
$\chi(y)\stackrel{{\rm def}}{=}0,$ for  $y\notin [0,B],$ and let
$$ \overline\chi(y)=\inf\,\{\,t>\chi(y):\;y+\xi(t)\in
[0,B]\,\},\qquad \overline X(y)=y+\xi(\overline\chi(y))\in
[0,B],\quad y\in \mathbb{R}
$$
be,respectively,  the  first entry time  into the interval $ [0,B] $
by the process $ y+\xi(t) $ and the value of the process at this
time. The integral transforms of the joint distribution of
$\{\overline\chi(y),\overline X(y)\},$ $y\in \mathbb{R} $ for $s>0$
satisfy the following equations
\begin{align}
b^v(du,s)&=
 e^{-vc(s)}\left(1-\frac{c(s)}{\lambda}\right)T(s)^{-1}P(\lambda,du),
\qquad  v>0,\notag\\
b_v(du,s)&= m_{v}^s(du)+
e^{Bc(s)}\left(1-\frac{c(s)}{\lambda}\right) \hat
T_v^{\,s}(c(s))\,T(s)^{-1}P(\lambda,du),
\qquad  v>0,\notag\\
b(y,du,s)&= \left(1-\frac{c(s)}{\lambda}\right)T(s)^{-1} \left      
[e^{c(s)(B-y)}-\frac{R_{B-y}(s)}{\hat R_B(\lambda,s)}\,
\frac{e^{-B(\lambda-c(s))}}{\lambda-c(s)}\right]
P(\lambda,du)\notag\\
&+\frac1{\lambda}\,\frac{R_{B-y}(s)}{\hat R_B(\lambda,s)}
(T(s)^{-1}-1)\,P(\lambda,du), \qquad  y\in [0,B],
\end{align}
where
\begin{align*}
& m_x^s(du)= E\,[\,e^{-s\tau^x};\;T^x\in\,du\,],\quad \hat
T_x^s(c(s))= E\,[\,e^{-s\tau^x-c(s)T^x};\;T^x>B\,],
\qquad x\ge0,\\
&\hat T_{\gamma}^s(c(s))= \lambda\int_0^{\infty}e^{-\lambda x}\hat
T_{x}^s(c(s))\,dx,\qquad
T(s)=1-\left(1-\frac{c(s)}{\lambda}\right)\hat T_{\gamma}^s(c(s))
e^{-B(\lambda-c(s))}.
\end{align*}
\end{theorem}
\begin{proof} \rm
We apply now the equalities (17) of Theorem 2 to obtain the formulae
(23). For  this we have to calculate  for the Poisson process with a
negative component  the kernels   $ Q_{\pm}(v,dl,s),$ and the
successive iterations  $ Q_{\pm}^{\,(n)}(v,dl,s),$ $n\in\mathbb{N},$
and the series   $ Q_{\pm}^{\,s}(v,dl).$ Utilizing the defining
formula of the kernels   (21) and the formulae (10) yields
\begin{align}
&Q_{+}(v,dl,s) = e^{-vc(s)}\left(1-\frac{c(s)}{\lambda}\right)
e^{-\lambda B}
 E[e^{-s\tau^{\gamma}};T^{\gamma}-B\in dl],                           
\qquad & v>0,\notag\\
&Q_{-}(v,dl,s) = \hat
T_v^s(c(s))\,e^{-B(\lambda-c(s))}\left(1-\frac{c(s)}{\lambda}\right)
\lambda\,e^{-\lambda l}\,dl,\qquad & v>0,
\end{align}
where $\hat T_x^s(c(s))= E[e^{-s\tau^x-c(s)T^x};T^x>B\,],$ $ x\ge0.$
Using the defining formula (19) for the successive iterations and
the method of mathematical induction it follows from  (24) that for
$n\in\mathbb{N},$
\begin{align*}
&Q_{+}^{\,(n)}(v,dl,s) =
e^{-vc(s)}\left(1-\frac{c(s)}{\lambda}\right) e^{-\lambda B}
\left(\tilde T_{\gamma}^s(c(s))\right)^{n-1}
E[e^{-s\tau^{\gamma}};T^{\gamma}-B\in dl],
,\notag\\
&Q_{-}^{\,(n)}(v,dl,s) = \hat
T_v^s(c(s))\,e^{-B(\lambda-c(s))}\left(1-\frac{c(s)}{\lambda}\right)
\left(\tilde T_{\gamma}^s(c(s))\right)^{n-1}\lambda\,e^{-\lambda
l}\,dl,
\end{align*}
where
\begin{align*}
\tilde
T_{\gamma}^s(c(s))=e^{-B(\lambda-c(s))}\left(\lambda-c(s)\right)
\int_0^{\infty}e^{-\lambda x}\hat T_x^s(c(s))\,dx.
\end{align*}
The series   $ Q_{\pm}^{\,s}(v,dl)$  of the successive iterations
$Q_{\pm}^{\,(n)}(v,dl,s)$ (see (18)) are just the geometric series,
and their sums are given by
\begin{align*}
& Q_{+}^{\,s}(v,dl) = \delta(v-l)\,dl+
e^{-vc(s)}\left(1-\frac{c(s)}{\lambda}\right) e^{-\lambda B}
 T(s)^{-1} E\,[\,e^{-s\tau^{\gamma}};\;T^{\gamma}-B\in dl\,],
\quad & v>0,\notag\\
& Q_{-}^{\,s}(v,dl) =\delta(v-l)\,dl+ \hat
T_v^s(c(s))\,e^{-B(\lambda-c(s))}\left(1-\frac{c(s)}{\lambda}\right)
T(s)^{-1}\lambda\,e^{-\lambda l}\,dl, \quad & v>0,
\end{align*}
where $ T(s)=1- \tilde T_{\gamma}^s(c(s)).$ Substituting   in the
equalities   (17) of Theorem 2  the expressions for the functions $
 Q_{\pm}^{\,s}(v,dl),$ and the expressions for the functions  $
 E\,[\,e^{-s\chi(y)};\;X(y)\in dv,\,A^{\,B}],$ $
E\,[\,e^{-s\chi(y)};\;X(y)\in dv,\,A_{0}],$  which are given by the
formulae of Theorem 1, we obtain the formulae (23) of  Theorem 3.
\end{proof}
\noindent{\bf Acknowledgements.}  The  authors  acknowledge support
of  the Belgian Federal Science Policy. (Interuniversity Attraction
Pole Programme P6/STADEC).


\begin{thebibliography}{9}

\bibitem{Bert1996} %
 Bertoin,  J. (1996). {\itshape L\'{e}vy processes,} Cambridge
University Press.

\bibitem{Bert1997}
 Bertoin,  J. (1997). Exponential decay  and ergodicity  of
 completely  assymetric  L\'{e}vy process in a finite interval. {\itshape Ann.Appl.
 Probab,} {\bfseries 7,} 156-169.

\bibitem{Br1}
Borovkov,  A.A. (1976).  {\itshape  Stochastic  processes in
queueing theory, }Springer Verlag.


\bibitem{DiP1}%
Ditkin,  V.A., Prudnikov,  A.P. (1966). {\itshape  Operational
calculus,} Moscow,  Russian edition.


\bibitem {Emery}
Emery, D.J. (1973).  Exit problem for a spectrally positive process.
{\it Adv. Appl. Prob.} \: {\bfseries 5, } 498--520.

\bibitem{Sk2}%
 Gihman, I.I.,  Skorokhod, A. V. (1975). {\itshape Theory of stochastic
processes, } Vol. 2, Springer Verlag, translated  from the Russian
by S. Kotz.
\bibitem{BGu} %
 Bratiychuk, N.S., Gusak D.V. (1990). {\itshape Boundary problems
 for processes with independent increments,} Kiev, Naukova Dumka, in Russian.



\bibitem{2Ka2}  Kadankov, V. F., Kadankova, T. V. (2005).
On the distribution of the first exit time from  an interval and the
value of overshoot  through the boudaries   for  processes with
independent increments and random walks. {\it Ukr. Math. J.}
{\bfseries 10 (57),} 1359-1384.


\bibitem{Kad6} %
Kadankov, V.F.   Kadankova, T.V. (2004). On the disribution of
duration of stay in an interval of the semi-continuous process with
independent increments. {\it Random Oper. and Stoch. Equa. (ROSE)}
{\bfseries 12(4),} 365--388.

\bibitem{Kad3}%
Kadankova, T.V. (2003). On the distribution of the number of the
intersections of a fixed interval by the semi-continuous process
with independent increments. {\it Theor. of Stoch. Proc.}\:
{\bfseries 1-2,} 73--81.


\bibitem{Kad2}
Kadankova, T.V. (2004). On the joint distribution of supremum,
infimum and  the magnitude of a process with independent increments.
 {\it Theor. Prob.  and Math.  Statist.}    {\bfseries 70,} 54--62.



\bibitem{Kypr2003}
Kyprianou, A.E. (2003).   A martingale review of some fluctuation
theory for spectrally negative L\'{e}vy processes, Research report,
Utrecht University.


\bibitem{PeRog} %
Pecherskii, E.A.,  Rogozin B.A. (1969).  On joint distributions of
random variables assosiated with fluctuations of a process with
independent increments. {\it Theor. Prob. and its Appl. }{\bfseries
14,} 410--423.


\bibitem{Pet1}
Petrovskii, I.G. (1965). {\itshape  Lectures on the theory of
integral equations,} Moscow, Nauka, Russian edition.

\bibitem{Pistor2004}
 Pistorius,  M.R. (2004).   A potential theoretical review of some
exit problems  of spectrally negative L\'{e}vy processes. {\it
S\'{e}minaire de Probabilit\'{e}s}, {\bf 38,} 30--41.

\bibitem{Pist2004}
Pistorius, M. R.  (2004). On exit and ergodicity of the spectrally
negative L\'{e}vy process reflected in its infimum.
 {\it J. Theor. Prob.} \:  {\bfseries
17,} 183--220.

\bibitem{Rog1}%
Rogozin, B.A. (1966).  On distributions of functionals related to
boundary problems for processes with independent increments.  {\it
Theor. Prob.  and its Appl.} {\bfseries 11(4),} 656-670.


\bibitem{Su} %
Suprun, V.N. (1976). Ruin problem and the resolvent of a terminating
process with independent increments.  {\it Ukr. Math. J.} {\bfseries
28(1),} 53--61,  (English transl.).

\bibitem{SuSh}
 Suprun, V.N.,  Shurenkov, V.M. (1976).  On the resolvent of a
process with independent increments terminating at the moment when
it hits the negative real semiaxis. In {\itshape Studies in the
Theory of Stochastic processes,} Institute of Mathematics, Academy
of Sciences of UKrSSR, Kiev, 170-174.

\bibitem{Zol1}%
Zolotarev, V. M. (1964).  The first passage time of a level and the
behaviour at infinity for a class of processes with independent
increments. {\it  Theor. Prob. and its Appl.}  {\bfseries 9(4),}
653--664.

\end{thebibliography}
\end{document}